\documentclass[12 pt]{article}
\usepackage{amssymb}
\usepackage[centertags]{amsmath}
\usepackage{amscd, amsthm}

\newtheorem{fact}{fact}[section]
\newtheorem{esempio}[fact]{Example}
\newtheorem{thm}[fact]{Theorem}
\newtheorem{lemma}[fact]{Lemma}
\newtheorem{prop}[fact]{Proposition}
\newtheorem{corollario}[fact]{Corollary}
\newtheorem{defini}[fact]{Definition}
\newtheorem{osserva}[fact]{Remark}

\begin{document}
\title{A theorem of complete reducibility for exponential polynomials}
\author{P. D'Aquino\footnote{Department of Mathematics, Seconda Universit\`{a} di Napoli, Via Vivaldi
43, 81100 Caserta, paola.daquino@unina2.it} and G.
Terzo\footnote{Department of Mathematics, Seconda Universit\`{a}
di Napoli, Via Vivaldi 43, 81100 Caserta,
 giuseppina.terzo@unina2.it}}
\maketitle

\begin{abstract}
In this paper we give a factorization theorem for the ring of
exponential polynomials in many variables over an algebraically closed field of characteristic 0
with an exponentiation.
This is a generalization of the factorization theorem due to Ritt
in \cite{ritt}.

\end{abstract}

\section{Introduction}
 Ritt  in \cite{ritt}  was the first to
consider a factorization theory for exponential polynomials over the complex field of the form
\begin{equation}
\label{Rittpoly}
f(z) = a_1e^{\alpha_1z} + \ldots +a_ne^{\alpha_nz},
\end{equation}
where
$a_i, \alpha_i \in \Bbb C.$ Contributions by  Gourin and  Macoll
in \cite{gourin} and \cite{macoll} gave refinements of Ritt's
result  for exponential polynomials of the form

\begin{equation}
\label{vdppoly}
f(z) = p_1(z)e^{\alpha_1z} + \ldots +p_n(z)e^{\alpha_nz},
\end{equation}
where $\alpha_i$ are complex numbers and  $p_i(z)\in \mathbb
C[z]$. Only in the mid 1990's  van der Poorten and Everest obtained a
factorization theorem  which applies to exponential
polynomials of the form  (\ref{vdppoly})
over any algebraically closed  field of characteristic
$0$. In \cite{vdpoorteneverest} they state (without proving)  that their result applies to a more general setting of a group ring
$R[G]$, where $R$ is a unique factorization domain and $G$ is a
divisible torsion-free ordered abelian group.

\par
The basic idea introduced
originally by Ritt is that of reducing the factorization of an
exponential polynomial to that of a classical polynomial in many
variables allowing fractional powers of the variables. In general, if we consider an
irreducible  polynomial $Q(y_1, \ldots, y_p)$ over a field it can  happen that for some positive integers $\mu_1, \ldots,
\mu_{p}$ the polynomial $Q(y_{1}^{\mu_1}, \ldots, y_p^{\mu_p})$
becomes reducible.
This may occur when we work with exponential
polynomials.  Ritt and Gourin saw the relevance, in terms of
factorization, of understanding the ways in which an
irreducible polynomial $Q(y_1, \ldots, y_p)$ becomes reducible
once the variables are replaced with their powers. The
determination of the integers $\mu$'s for which the reducibility
occurs is a crucial step in their results. More precisely, they
give a uniform bound for the number of irreducible factors of
$Q(y_{1}^{\mu_1}, \ldots, y_p^{\mu_p})$ which depends only on the
degree of $Q(y_1, \ldots, y_p)$. van der Poorten in
\cite{vanderpoorten} refined the bound by proving that it depends
only on the degrees of two of the $n$ variables. \par
In this paper, we generalize  the main ideas due to Ritt to the ring of
exponential polynomials in many variables over an
algebraically closed field of characteristic $0$ with an exponentation. We produce a  theorem of complete reducibility for  exponential polynomials with any iterations of
exponentiation.

\section{$E$-polynomial ring}
\label{anelloesponenziale}

\begin{defini} An exponential ring, or $E$-ring, is a pair $(R,
E)$ with $R$  a ring (commutative with 1) and $$E: (R, +)
\rightarrow ({\cal U}(R), \cdot)$$ a map of the additive group of
$R$ into the multiplicative group of units of $R$, satisfying
\begin{enumerate}
\item $E(x + y) = E(x) \cdot E(y)$ for all $x, y \in R$ \item
$E(0) = 1.$
\end{enumerate}
\end{defini}

$(K, E)$ is an $E$-field if $K$ is a field.\\

\par We now recall the construction of the ring of exponential polynomials ($E$-polynomials) over an $E$-field $(K, E)$  and some
of its basic properties. The construction  is well known (see  \cite{van} or \cite{mac}),
but for the proof  of the main result of this paper we prefer to go through the details. From now on we will be working with $(K,E)$ an algebraically closed field of characteristic 0 with an exponentiation, unless otherwise specified.

\par The ring of exponential polynomials ($E$-polynomials) in
the indeterminates $\overline{x} = x_1, \ldots, x_n$ is an
$E$-ring constructed as follows by recursion. We construct three
sequences:
\begin{enumerate} \item $(R_{k}, +, \cdot)_{k\geq-1}$ are rings; \item
$(B_k, +)_{k\geq0}$ are torsion free divisible abelian groups;
\item $(E_k)_{k\geq-1}$ are partial
$E$-morphisms.\\
\end{enumerate}
\textbf{Step 0:} Let\\
 $R_{-1} = K;$\\
 $R_0 =(K[\overline x], +, \cdot);$\\
 $B_0=\langle \overline x \rangle$,  the ideal generated by $\overline x$. So $R_0 = R_{-1} \oplus B_0;$ \\
  $E_{-1}: R_{-1}
\longrightarrow
 R_0,$ is the composition of the initial $E$-morphism over $K$ with the immersion of $K$ into
$K[\overline x].$\\
\textbf{Inductive step:}\\
Suppose $k \geq 0$ and $R_{k-1}$, $R_{k}$, $B_k$ and $E_{k-1}$ have been defined in such a way that $R_k$ as additive group is
 $$ R_k = R_{k-1} \oplus B_{k}, \mbox{  and  } E_{k-1}: (R_{k-1}, +) \rightarrow ({\cal
U}(R_k),
 \cdot),$$
where ${\cal U}(R_{k})$ denotes the set of units in $R_{k}.$
\par Let
\[ t: (B_k, +) \rightarrow (t^{B_k}, \cdot) \]
be a formal isomorphism. This is used in order to convert an additive group into a multiplicative copy of it. Define \\
\begin{center}
$R_{k+1} = R_{k} [t^{B_{k}}]$ (as group ring over $R_k$).\\
\end{center}
Therefore \[ R_{k} \mbox{ is a subring of } R_{k+1}
\] and as additive group
$$R_{k+1} = R_{k} \oplus B_{k+1},$$
where $B_{k+1}$ is the $R_k$-submodule of $R_{k+1}$ freely
generated by $t^{b},$ with $b \in B_{k}$ and $b \neq 0$ (this last
condition ensures that $B_{k+1}$ does not coincide with
$R_{k+1}$).
\par We define
$$E_{k}: (R_{k}, +) \rightarrow ({\cal U}(R_{k+1}), \cdot)$$
as follows  $E_k(q) = E_{k-1}(r) \cdot t^{b},$ for $q = r + b,$
$r \in
R_{k-1}$ and $b \in B_k$.\\
\par We have constructed a chain of partial $E$-rings (the domain of
exponentiation of $R_{k+1}$ is $R_k$)
 \[R_0 \subset R_1 \subset R_2 \cdots \subset R_k \subset
\cdots\]
\noindent Then the $E$-polynomial ring is\\
$$K[\overline{x}]^E = \lim_{k} R_k = \bigcup_{k=0}^{\infty}R_k$$
and the $E$-ring morphism defined on $K[\overline{x}]^E$ is defined
as follows $$ E(q) = E_k(q),  \mbox{ if } q \in R_k \mbox{ and } k \in
\Bbb N.$$ \par Notice that each $R_{k+1}$ as additive group is the
direct sum $K \oplus B_0 \oplus B_1 \oplus \ldots \oplus B_{k+1}.$
Moreover, as an additive group $K[\overline{x}]^E$ is $$K \oplus B_0 \oplus B_1 \oplus \ldots \oplus
B_{k+1} \oplus \ldots.$$ \par
For all $k$ the group
ring $R_{k+1}$ can be viewed in the following different ways
$$R_{k+1} \cong R_0[t^{B_0 \oplus B_1 \oplus \ldots \oplus B_k}];$$
$$R_{k+1} \cong R_1[t^{B_1 \oplus \ldots \oplus B_k}];$$
$$...$$
$$...$$
$$R_{k+1} \cong R_k[t^{B_k}].$$ Moreover, $K[\overline x]^E =
R_{0}[t^{B_0 \oplus B_1 \oplus \ldots \oplus B_k \ldots}],$ i.e.
$K[\overline x]^E$ is a group ring constructed  over a
unique factorization domain, $R=K[\overline x](=R_0)$ and a torsion free divisible
abelian group $G=t^{B_0 \oplus B_1 \oplus \ldots \oplus B_k
\ldots}$ (a $\mathbb Q$-vector space).

If  $f \in K[\overline x]^E$ then there is a unique $k$ such that $f\in R_{k+1} \cong R_{0}[t^{B_0 \oplus B_1 \oplus \ldots \oplus B_k }]$ and $f\not\in R_k$.
In this case $f$ can be written uniquely as
\begin{equation}\label{polinomio} f= \sum_{h=1}^{N} a_ht^{\alpha_h},
\end{equation}
where $a_{h} \in R_{0}$ and the $\alpha_{h} \in B_1\oplus
\ldots \oplus B_k.$

\begin{defini}\par Let $f \in K[\overline x]^E$ be as in (\ref{polinomio}).
 The support of $f,$ denoted by $supp(f),$ is the $\mathbb Q$-vector space generated by $\alpha_{1}, \ldots, \alpha_{N}.$
 \end{defini}

\begin{osserva}\label{ordine} \rm{From the construction $B_1\oplus \ldots \oplus B_k$ is a torsion
free abelian group, and using a compactness argument it can be made
into an ordered group with order $<.$ Without loss of generality, at the cost of the
renumbering the $\alpha'$s we can assume that
$$\alpha_{1} < \alpha_{2} < \ldots < \alpha_{N}.$$}\end{osserva}


We recall some properties which are preserved from the
$E$-ring $(K,E)$ to $K[\overline x]^E$, the ring of $E$-polynomials.

\begin{prop} If $K$ is an integral $E$-domain of characteristic 0 then
the $E$-polynomial ring $K[\overline x]^E$ is an integral domain
whose units are of the form $uE(\alpha),$ where $u$ is a unit of
$K,$ and $\alpha$ is an exponential polynomial.
\end{prop}

Note that if the exponential map is surjective (as over $\Bbb C$) then the units correspond to purely exponential terms.



\begin{defini}
An exponential polynomial $f \in K[\overline x]^E$ is
irreducible if there are no non-units g and h with $f = gh.$
\end{defini}

\section{Associate polynomial}\label{associato}

Using the ideas of Ritt we associate
to any exponential polynomial a classical polynomial in many
variables over a unique factorization domain. This is a crucial step in order to obtain a
factorization theorem for the exponential  ring $K[\overline x]^E$.
\par
Let $f\in K[\overline x]^E$
\begin{equation}
f = \sum_{h=1}^{N} a_ht^{\alpha_h},
\end{equation}
where $a_{h} \in K[\overline x]$ and $\alpha_{h} \in B_1\oplus
\ldots \oplus B_k.$
If
 $\{ \mu_1,\dots ,\mu_p\} $  is  a $\Bbb Q$-basis
of $supp(f)$ then
$$\alpha_h= r_{h1}\mu_1+\ldots +r_{hp}\mu_p$$ for all $h=1,\ldots
, N$, and $r_{hj} \in \Bbb Q.$ Notice  that we can always find a
new basis for $supp(f)$ so that $r_{hj}$'s are integers. Indeed, if
$M$ is the least
 common multiple of the denominators of $r_{hj}$'s then there are integers $s_{i1},\ldots ,s_{ip}$
 such that $$\alpha_h= s_{h1}\frac{\mu_1}{M}+\ldots + s_{hp}\frac{\mu_p}{M}$$
 for $h=1,\ldots ,N$.
 Let $\nu_{1}=\frac{\mu_1}{M}, \ldots, \nu_{p}=\frac{\mu_p}{M}$. Clearly, $\nu_{1}, \ldots ,\nu_{p}$
 is a new basis  of $supp(f)$, and each $\alpha_h$ can be expressed as a linear combination of
 $\nu_{1}, \ldots ,\nu_{p}$ with integer coefficients.

Let $y_1=t^{\nu_1},\ldots ,y_p=t^{\nu_p}.$ Notice that the linear independence of  $\nu_i$'s implies the algebraic independence  of $t^{\nu_1},\ldots, t^{\nu_p}.$ By expressing each $\alpha_i$ in terms of the $\nu_j$'s we have that $f$ is transformed into a classical Laurent polynomial $Q$ over $K[\overline x]$ in the variable $y_1, \ldots, y_p.$ This means that  we can write $Q$ as a product of a polynomial in the $y_i'$s and a quotient of a  monomial in the $y_i'$s.

We have the following correspondence:
$$f \in K[\overline x]^E \rightsquigarrow Q(y_1, \ldots, y_p) \in R[y_1, \ldots, y_p]$$
where $R=K[\overline x]$. We will refer to $Q(y_1, \ldots, y_p)$ as the associate polynomial of $f$.


\begin{osserva}\rm{The correspondence between $f$ and $Q$ holds modulo a monomial in
$y_1, \ldots, y_p$,  which corresponds to an invertible element of
$K[\overline x]^E $, and it does not have any consequence on the
factorization of $f$.}
\end{osserva}

\begin{defini}
An exponential polynomial $f$ is  simple if $\dim(supp(f)) = 1.$
\end{defini}

An example of a simple polynomial in $\Bbb C[x]^E$ is $$ sin(2\pi x) =
\frac{e^{2\pi ix}-e^{-2\pi i x}}{2i}$$

\medskip

A classical polynomial $Q(y_1, \ldots, y_p) $ is {\it essentially
$1$-variable} if there are monomials $\tau_1, \tau_2$ in $y_1,
\ldots, y_p$ such that $Q(y_1, \ldots, y_p) =\tau_1P(\tau_2)$,
where $P$ is a polynomial in just one variable. Hence, a simple
exponential polynomial has associated an essentially $1$-variable
polynomial. In other words, a simple polynomial is a polynomial in
$e^{\mu x}$, where $\mu$ is a generator of  $supp(f).$

\medskip

An example of an essentially $1$-variable polynomial is $$Q(x,y)= x^2y(3x^3y^9-2x^2y^6+1)= \tau_1P(\tau_2),$$ where $\tau_1=x^2y$, $\tau_2=xy^3$ and   $P(z)=3z^3-2z^2+1.$

\subsection{Classical polynomials in many variables}

The correspondence between exponential polynomials and
classical polynomials implies  that there are connections
between factorizations of one in terms of  factorizations of the
other one, and vice versa. 
Ritt in \cite{ritt} detected  that the classical factorization of
a polynomial in many variables could not exhaust  all the
factorizations of the exponential polynomial, and he introduced
factorizations in fractional powers of the variables. In the
context of exponential polynomials over an exponential field
fractional powers of an exponential do make sense.  Refinements of
his ideas where obtained by Gourin in \cite{gourin},  and more
recently by van der Porten and Everest in \cite{vanderpoorten} and
\cite{vdpoorteneverest}.

\smallskip

\par
Let  $Q(y_1, \ldots, y_{p}) \in F[y_1, \ldots, y_{p}]$ be an irreducible polynomial over $F$, where $F$ is a field or a unique factorization domain. It can happen that for some ${\mu_1}, \ldots, {\mu_p} \in
\Bbb N-{0},$ $Q(y_1^{\mu_1}, \ldots, y_p^{\mu_p})$ becomes
reducible. For example, if $Q(x,y)=x-y$ then $Q(x^3,y^6)=(x-y^2)(x^2+xy^2+y^4)$.
\par

\begin{defini} A polynomial $Q(y_1, \ldots, y_{p})$ is power irreducible
(over $F$) if for each ${\mu_1}, \ldots, {\mu_p} \in
\Bbb N-{0},$ $Q(y_1^{\mu_1}, \ldots, y_p^{\mu_p})$ is irreducible.
\end{defini}

\begin{defini}
A polynomial $Q(y_1, \ldots, y_{p})$ is primary in the variable $y_{i}$ if
$$Q(y_1, \ldots, y_{p}) = P(y_{1}, \ldots,
y_{i}^d, \ldots, y_p)$$ for some polynomial $P(y_1, \ldots, y_{p}) \in K[y_1, \ldots, y_{p}]$  implies $d =1.$ Equivalently,
$Q(y_1, \ldots, y_{p})$ is primary in $y_{i}$ if the g.c.d.
 of the exponents of $y_i$ in all terms of $Q$ is $1$.
\end{defini}

\begin{defini}
A polynomial $Q(y_1, \ldots, y_{p})$ is primary if it is primary in each
variable.
\end{defini}

\begin{esempio}
$Q(x,y) = 3x^2y - 5y^3 + x^3$ is primary in both $x$
and $y.$ On the contrary, $R(x, y) = 3x^2y - 5y^3 + x^4$ is not
primary in $x$ since $R(x, y) = P(x^2, y)$, where
$P(x, y) = 3xy - 5y^3 + x^2$
\end{esempio}

\begin{osserva}
{\rm Notice that if $Q(y_1, \ldots, y_{p})$ is a non primary
polynomial then  there exists a unique $p$-tupla $t_1, \ldots,
t_p$ of positive integers such that $$Q(y_1, \ldots, y_{p})=
P(y_1^{t_1}, \ldots, y_p^{t_p})$$ where $P(y_1, \ldots, y_{p})$ is
primary. We will see that it is not restrictive to work with
primary polynomials in connection to the factorization of an
exponential polynomial.}
\end{osserva}

\section{Factorization Theorem}

In this section we analyze factorizations of $f$ in terms of
factorizations of  $Q(y_1, \ldots, y_p)$ in fractional powers of
$y_1, \ldots, y_p$ over a unique factorization domain $R$ of
characteristic $0$ containing all roots of unity.



\par

\medskip

If fractional powers are permitted then a binomial as $y - 1$
defined over an algebraically closed field $K$ may have infinitely
many factors. Indeed,  $y^{\frac{1}{k}} - \epsilon$, where
$\epsilon$ is a $k$th root of unity, is a factor of $y - 1$ for
any positive integer $k.$ So, any general discussion of
factorization in fractional power must avoid such polynomials
because of the infinitely many factors they may have.

We recall that a simple polynomial $f$ has associated an
essentially $1$-variable polynomial $Q(y)$ which factorizes over
the algebraic closure $R^{alg}$ of $R(=K[\overline{x}])$ into a
finite number of polynomials of the form $1 + a y$ with $a \in
R^{alg}.$ As observed before if fractional powers of the variables
are allowed, $1 + ay$ has factors of the form $1
+a'y^{\frac{1}{k}}$ for $k= 1, 2, 3...$, and $a'\in R^{alg}.$ For
this reason in the factorization theorem fractional powers of the
variable will be avoided  for simple polynomials.

\subsection{Main result}

\par Ritt and Gourin (1927-1930), and van der Poorten
(1995) studied factorizations of $Q(y_1, \ldots, y_p)$ into
primary irreducible polynomials in the ring generated over an
algebraically closed field $K$ of characteristic $0$ by all
fractional powers of $y_1, \ldots, y_p.$ More precisely they
proved the following theorem:

\begin{thm}[\cite{ritt},\cite{gourin}]
\label{finitenessGourin}
There is a uniform bound for the number of primary,
irreducible factors of
$$Q(y_1^{\mu_1},\ldots ,y_p^{\mu_p})$$ for $Q(y_1, \ldots, y_p)$
not effectively 1-variable, and arbitrary $\mu_1,\ldots ,\mu_p \in
\Bbb N_{+}.$ The bound depends only on
$$M = max\{ d_{y_1}, \ldots, d_{y_p}\}$$
where $d_{y_i}$ denotes the $y_i$-degree.
\end{thm}

\noindent van der Poorten in \cite{vanderpoorten} obtaines a  bound in terms of the  degrees of only two variables.
\par
As a consequence the following corollary holds.

\begin{corollario} The number of irreducible factors of $Q(y_1, \ldots, y_p)$ in
fractional powers of $y_1, \ldots, y_p$ is finite.\end{corollario}

The factorization of an exponential polynomial $f$ is obtained via
the factorization of the associate polynomial $Q$ in fractional
powers of the variables. Indeed, finite factorizations of $f$
generate factorizations of $Q$ in fractional powers of the
variables, and vice versa.  One implication depends on the
following main lemma.


\begin{lemma}
Let $f(\overline x) \in K[\overline x]^E$ and suppose that
$f(\overline x) = g(\overline x) \cdot h(\overline x),$ where
$g(\overline x), h(\overline x) \in K[\overline x]^E.$ Then
$supp(g), supp(h)$ are contained in $supp(f).$
\end{lemma}

\textbf{Proof:} Let $f(\overline x)= a_{1}t^{\alpha_1} +
a_{2}t^{\alpha_2} + \ldots + a_{N}t^{\alpha_N}$, and suppose that

\begin{equation}\label{divisible} f =
(b_{1}t^{\beta_1} + b_{2}t^{\beta_2} + \ldots + b_{M}t^{\beta_M})
\cdot (c_{1}t^{\gamma_1} + c_{2}t^{\gamma_2} + \ldots +
c_{S}t^{\gamma_S}),\end{equation}

where $g = b_{1}t^{\beta_1} +
b_{2}t^{\beta_2} + \ldots + b_{M}t^{\beta_M}$ and $h =
c_{1}t^{\gamma_1} + c_{2}t^{\gamma_2} + \ldots +
c_{S}t^{\gamma_S}.$ \par Without loss of generality  (see Remark \ref{ordine}) we can assume
$$\alpha_{1} < \alpha_2 < \ldots < \alpha_N,$$
$$\beta_{1} < \beta_2 < \ldots < \beta_M,$$
$$\gamma_{1} < \gamma_2 < \ldots < \gamma_S.$$
First of all we show that $supp(h)  $ is contained in the $\mathbb
Q$-space generated by $supp(f)$ and $supp(g)$, i.e. $supp(h)
\subseteq \langle supp(f), supp(g)\rangle_{\Bbb Q}.$ Suppose by
contradiction that $supp(h) \nsubseteq \langle supp(f), supp(g)
\rangle_{\Bbb Q}.$ Let $\gamma$ be the maximum of $\gamma_{1},
\gamma_2, \ldots, \gamma_S,$  such that $\gamma$ is not in
$\langle supp(f), supp(g) \rangle_{\Bbb Q}.$ If there are no
$\beta_{i}$ in $g$ and $\gamma_{j}$ in $h$ such that $\beta_{i} +
\gamma_{j} = \beta_M + \gamma,$ this means that the term
$t^{\beta_M + \gamma}$ does not cancel out and $\beta_M + \gamma =
\alpha_{k}$ for some $k =1, \ldots, N,$ so $\gamma \in \langle
supp(f), supp(g) \rangle_{\Bbb Q}.$ Otherwise, suppose that there
are some $\beta_{i}$ and $\gamma_j$ such that $\beta_M + \gamma =
\beta_{i} + \gamma_{j}.$ Since $\beta_M > \beta_{i},$ then it must
be $\gamma_i > \gamma$, and so $\gamma_i \in \langle supp(f),
supp(g) \rangle_{\Bbb Q}.$ Hence, $\gamma = -\beta - \beta_{i} +
\gamma_i$ is in $\langle supp(f), supp(g) \rangle_{\Bbb Q}.$ In
both cases we have a contradiction. Hence,  $supp(h) \subseteq
\langle supp(f), supp(g) \rangle_{\Bbb Q}.$ We can repeat a
similar argument for $g$, and we obtain that $supp(g) \subseteq
\langle supp(f), supp(h) \rangle_{\Bbb Q}.$ The previous two
inclusions, clearly, imply that $$\langle supp(f), supp(h)
\rangle_{\Bbb Q} = \langle supp(f), supp(g) \rangle_{\Bbb Q} =
\langle supp(f), supp(g), supp(h) \rangle_{\Bbb Q}.$$ It remains
to prove that $$supp(g),supp(h) \subseteq supp(f).$$ Suppose by
contradiction that it is not. Let $d_1, \ldots, d_t$ be a $\Bbb
Q$-basis of $\langle supp(f), \break supp(g), supp(h)
\rangle_{\Bbb Q}.$ We can write $\beta_M, \gamma_s$ in terms of
this basis and in particular we can write their sum as a $\Bbb
Q$-linear combination of $d_1, \ldots, d_t$. The order of the
$\beta_i$'s and $\gamma_j$'s implies that the term $t^{\beta_M +
\gamma_S}$ does not cancel out and  so,  by relation
(\ref{divisible}), $\beta_M + \gamma_S = \alpha_N.$ This
contradicts the $\Bbb Q$-linear independence of $d_1, \ldots,
d_t.$ This completes the proof. \qed

\smallskip

Previous lemma implies the following corollary.
\begin{corollario} Suppose $f$  factorizes as $f=gh$.
Let $Q$, $P$ and $R$ be the associate polynomials to $f$, $g$ and
$h$, respectively. There is a monomial $\tau$ such that $Q=P'R'$
where $P'=P\tau$, $R'=\tau^{-1}R$.
\end{corollario}

\begin{osserva}
{\rm  If $f$ is a simple polynomial and $g$ divides $f$ then $g$ is also simple.

}
\end{osserva}

Factorizations of the associate polynomial  do not exhaust all the factorizations of $f$.
We need the following crucial result which is a generalization of Theorem \ref{finitenessGourin} for polynomials over a unique factorization domain contaning all roots of unity.  This ensures that the number of irreducible factors in fractional powers of the variables of a polynomial is finite.

\begin{thm}
\label{finiteness}
Let $Q(y_1,\ldots ,y_p)$ be a primary irreducible polynomial over $R$, a unique factorization domain of characteristic $0$ contaning all roots of unity. Assume also that $Q$ is  not essentially a
$1$-variable polynomial. Then $Q(y_1,\ldots ,y_p)$ has a factorization into primary irreducible
polynomials in the ring generated over $R$ by all fractional powers of  $y_1,\ldots ,y_p$, and
the number of these irreducible factors is finite.
\end{thm}

\medskip
We will not give the detailed proof of Theorem \ref{finiteness} since it can be obtained step by step following the corresponding proof  of the analogue result  in \cite{gourin}.
As already remarked the main difference of our results is that we work over a unique factorization domain $R(=K[\overline{x}])$ containing all roots of unity since $K$ is an  algebraically closed field while in \cite{gourin} the polynomials are over $\mathbb C$. In \cite{gourin}  a non essentially $1$-variable polynomial corresponds to a polynomial with at least three terms since the polynomials are over $\mathbb C$. In some of our arguments we will tacitly  work over the algebraic closure of $R$, and this is not a restriction with respect to the finiteness of the number of irreducible factors of a polynomial.

\smallskip

First of all notice that any factorization of $Q(y_1^{t_1}, \ldots ,y_p^{t_p})$
for $t_1, \ldots , t_p\in \mathbb N$ gives a factorization in fractional powers of the variables of $Q(y_1,\ldots ,y_p)$.
Hence the factorizations of $Q(y_1^{t_1}, \ldots ,y_p^{t_p})$ are relevant.
\par

\medskip

The main steps of the proof can be summarized as follows.

{\bf Step 1.}  For $i=1,\dots ,p$ let $\epsilon_i $ be a primitive $t_i$th root of unity, and consider the transformations $\tau_{\epsilon_i}: y_i\mapsto \epsilon_i^k y_i$ for $0\leq k<t_i$.
Let $G$ be the group generated by  $\tau_{\epsilon_i}$'s, for $i=1,\dots ,p$ and $0\leq k<t_i$.

The polynomial $Q(y_1^{t_1}, \ldots ,y_p^{t_p})$ is left unchanged by the action of $G$ on its irreducible factors.
This is proved by showing that from one irreducible factor of
$Q(y_1^{t_1}, \ldots ,y_p^{t_p})$  all the others can be obtained via the group $G$ of transformations
generated by $\tau_{\epsilon_i} $'s.

We recall that $Q$ is an irreducible polynomial and consists of more than two terms, since it is not essentially $1$-variable.
\medskip

{\bf Step 2.} If the irreducible factor $Q_1$ of $Q^{(t)}$ is primary then $Q_1$ also consists of more than two terms.

\medskip
{\bf Step 3.} If an irreducible factor $Q_1$ of
$Q^{(t)}=Q(y_1^{t_1}, \ldots ,y_p^{t_p})$ is primary then each
$t_j$ satisfies the relation $t_j\leq M^2$, where
$M=max(d_{y_1},\ldots ,d_{y_p})$.

\smallskip

This implies that there are only finitely many sets  of positive
integers $$t_{11},\ldots ,t_{1p};  t_{21},\ldots ,t_{2p}; \ldots
\ldots ;t_{n1},\ldots ,t_{np}$$ such that $Q(y_1^{t_{i1}},\ldots
,y_p^{t_{ip}})$, for $i=1, \ldots , n$ are reducible.

\begin{osserva}
\rm{ The above results have been obtained for $Q$ a primary polynomial. Suppose $Q(y_1,\ldots ,y_p)=P(y_1^{d_1},\ldots ,y_p^{d_p})$ where $d_j$ are positive integers and at least one is strictly greater than $1$, and $P(y_1,\ldots ,y_p)$ is primary. If $$t_{11},\ldots ,t_{1p};  t_{21},\ldots ,t_{2p}; \ldots \ldots ;t_{n1},\ldots ,t_{np}$$ are the only sets such that $P(y_1^{t_{i1}},\ldots
,y_p^{t_{ip}})$, for $i=1, \ldots , n$ are reducible then $$\tau_{11},\ldots ,\tau_{1p};  \tau_{21},\ldots ,\tau_{2p}; \ldots \ldots ;\tau_{n1},\ldots ,\tau_{np}$$  where $\tau_{ij}= \frac{t_{ij}}{m_{ij}}$ with $m_{ij}=gcd(t_{ij}, d_j)$,
are the only sets of positive integers such that $Q(y_1^{\tau_{i1}},\ldots
,y_p^{\tau_{ip}})$, for $i=1, \ldots , n$ are reducible. }
\end{osserva}

We are now in a position to prove the following factorization theorem.

\begin{thm}
\label{main} An element $f\not= 0$ of $K[\overline x]^E,$ where
$K$ is an algebraically closed field of characteristic 0, factors
uniquely up to units and associates, as a finite product of
irreducibles of $K[\overline x]$, a finite product of irreducible
of $K[\overline x]^E$ whose support is of dimension bigger than
$1$,   and a finite product of elements  $g_j$ of $K[\overline
x]^E$, where $supp (g_{j_1})\not= supp (g_{j_2})$, for $j_1\not=
j_2$ and whose supports are of dimension $1$.
\end{thm}

\emph{Proof:} Let $f(\overline{x})\in K[\overline{x}]^E$, and
$Q(y_1,\ldots , y_p)$ the associate polynomial. We distinguish between the
irreducible factors of $Q$ which are essentially $1$-variable polynomials  and those which are not. The first ones correspond to the simple
factors of $f$, and we will multiply those which have the same support in order to have all the factors of $f$ of dimension $1$ of different support.

Now we consider the irreducible factors of $Q$ which are not essentially $1$-variable polynomials, and we study the factorizations of them in fractional powers of the variables. We will show how to get the corresponding irreducible factors of $f$ from them.
\par
Let $P(y_1,\ldots , y_p)$ be an irreducible factor of $Q$ not essentially $1$-variable. If $P(y_1,\ldots , y_p)=
V(y_1^{n_1},\ldots ,y_p^{n_p})$, for some $n_1, \ldots ,n_p\in
\mathbb N$ and $V(y_1,\ldots, y_p)$   primary then necessarily $V(y_1,\ldots
, y_p)$ is  irreducible otherwise $P(y_1,\ldots ,y_p)$ would
be reducible. Once we substitute each $y_i$ by $t^{n_i\nu_i}$ (where $\nu_1, \ldots ,\nu_p$ are defined as in Section 3) we get a
factor of $f$. Let $r_1, \ldots ,r_p$ be
positive integers such that $$V(y_1^{r_1}, \ldots, y_p^{r_p}) =
V_1 \cdot \ldots \cdot V_q$$ where $V_j$ are primary and
irreducible, and $q$ is the maximum number of
irreducible primary factors.  Theorem \ref{finiteness} guarantees that such $q$ exists.
Indeed, there are only  finitely many $p$-tuples of positive integers $n_1, \ldots ,n_p$ such that $V(y_1^{n_1}, \ldots, y_p^{n_p})$ is reducible. Among these take the one with the highest number of factors which is $q$.

\smallskip

\noindent {\it Claim.} The exponential polynomial obtained by replacing $y_i$ by $t^{n_i\nu_i/r_i}$ in $V_j$ for any $j=1,\dots, q$ is an irreducible  factor of $f$.

\smallskip

Suppose that for some $i_0$ this is not the case for $V_{i_0}$. Then for
$s_i=mn_i\nu_i/r_i$, where $m=lcm(r_1, \ldots , r_p)$ we have that
$V_{i_0}(y_1^{s_1}, \ldots, y_p^{s_p})$ is reducible. Hence,
\newline $V(y_1^{r_1}, \ldots, y_p^{r_p})$ has more than $q$ irreducible
factors. We can replace $r_is_i$ by a submultiple $z_i$ for
$i=1,\ldots, p$  in order to have a polynomial $V(y_1^{z_1},
\ldots, y_p^{z_p})$ with more than $q$ primary irreducible
factors. This is a contradiction with the maximality of $q$.

We have completed the proof of the existence of a factorization of an exponential  polynomial $f$ as a finite product of
irreducible not effectively $1$-variable polynomials and a finite product of simple polynomials.  Now we have to prove that such
factorization is unique. \par Suppose that $f(\overline x) \in
K[\overline x]^E$ has two different  factorizations as
$$f(\overline x) = g_1(\overline x) \cdot \ldots \cdot g_l(\overline x)$$
$$f(\overline x) = h_1(\overline x) \cdot \ldots \cdot h_s(\overline x).$$
It is enough to prove that if $g(\overline x)$
divides $h(\overline x) \cdot l(\overline x)$ in $K[\overline x]^E$  and
$g(\overline x)$ has no factor in common with $h(\overline x)$
 then $g(\overline x)$ divides $l(\overline x).$ Suppose

 \begin{equation}
 \label{equality}
 g(\overline x) \cdot s(\overline x) = h(\overline x) \cdot l(\overline x)
 \end{equation}
  for some $s(\overline x)\in K[\overline x]^E$, and $(g(\overline x), h(\overline x)) = 1.$ Let $G(\overline y), H(\overline y), L(\overline y), S(\overline y) $ be the associate polynomials to $g(\overline x),  h(\overline x), k(\overline x), s(\overline x) $, respectively. Clearly, $(G(\overline y), H(\overline y)) = 1,$ since any non trivial common factor of $G(\overline y)$ and  $H(\overline y)$ would give a non trivial common factor of $g(\overline x)$ and $h(\overline x)$.

  We saw  that any factorization of an exponential polynomial induces a factorization of the corresponding associate polynomial, hence  (\ref{equality}) implies the following relation over a unique factorization domain
\begin{equation}
\label{EQUALITY}
G(\overline y) \cdot S(\overline y) = H(\overline y) \cdot L(\overline y).
\end{equation}

Since $G$ has no common factors with $H$ then $G$ divides $L$.  This imples that the exponential polynomial $g(\overline x) $ divides $l(\overline x).$ So the uniqueness of the factorization of $f$ follows. \qed

\begin{osserva}
{\rm   The uniqueness of the factorization  has the following two important consequences:
\begin{enumerate}
\item
If  $f$ in  $K[\overline x]^E$ is irreducible and with support of dimension more than $1$ then $f$ is prime. For, if $f$ divides $gh$ then by the factorization  theorem $f$ must occur in the factorization of one of $g$ or $h$.

\item
In Section 3 in order to construct the associate polynomial to $f$ it was necessary to fix a basis of
$supp(f)$. Clearly, different bases determine different associate polynomials. From the unique factorization of $f$ it follows that they differ by a monomial.

\end{enumerate}
}

\end{osserva}

\medskip
\noindent {\bf Acknowledgements: }The authors would like to thank A. Macintyre for many helpful discussions and insights.

\end{document}